\input amstex 
\documentstyle{amsppt}
\input bull-ppt
\keyedby{bull431/lic}
\define\subd{_{\roman D}}
\define\otr{\operatorname{Tr}}
\define\spec{\operatorname{spec}}

\topmatter
\cvol{29}
\cvolyear{1993}
\cmonth{October}
\cyear{1993}
\cvolno{2}
\cpgs{250-255}
\title Trace Formulae and Inverse Spectral Theory\\
for Schr\"odinger Operators \endtitle
\author F. Gesztesy, H. Holden, B. Simon, and Z. 
Zhao\endauthor
\shorttitle{Trace Formulae}
\address \RM{(F. Gesztesy and Z. Zhao)} Department of 
Mathematics,
University of Missouri, Columbia, Missouri 65211 \endaddress
\ml\nofrills{\it E-mail address}, \RM{F. Gesztesy:}
mathfg\@mizzou1.missouri.edu
\mlpar
{\it E-mail address}, \RM{Z. Zhao:} 
mathzz\@mizzou1.missouri.edu
\endml
\address \RM{(H. Holden)} Department of Mathematical 
Sciences,
The Norwegian Institute of Technology, University of 
Trondheim, 
N-7034 Trondheim, Norway\endaddress
\ml holden\@imf.unit.no\endml
\address \RM{(B. Simon)} Division of Physics, Mathematics, 
and Astronomy, California Institute of Technology, 253--37,
Pasadena, California 91125\endaddress
\date February 2, 1993\enddate
\subjclass Primary 34A55, 34L40\endsubjclass
\thanks The first author
gratefully acknowledges support by the Norwegian Research
Council for Science and the Humanities (NAVF) and by 
Caltech.
The second author
is grateful to the Norwegian Research Council for Science 
and
the Humanities (NAVF) for support.
The research of the third author was partially supported 
by USNSF
Grant DMS-9101716\endthanks
\abstract We extend the well-known trace formula for 
Hill's equation
to general one-dimensional Schr\"odinger operators. The 
new function $\xi$,
which we introduce, is used to study absolutely continuous 
spectrum
and inverse problems.\endabstract
\endtopmatter

\document

In this note we will consider one-dimensional 
Schr\"odinger operators
$$
H=-\frac{d^2}{dx^2}+V(x)\quad\text{on}\ L^2(\Bbb R;dx)\tag 
1S
$$
and Jacobi matrices
$$
(hu)(n)=u(n+1)+u(n-1)+v(n)u(n)\quad\text{on}\ l^2(\Bbb 
Z).\tag 1J
$$
We will suppose that $V(x)$ is continuous and bounded 
below and
$v(n)$ is bounded.

In the analysis of the inverse problem for $H$ when $V$ is 
periodic
$(V(x+L)=V(x))$, a crucial role is played by a trace 
formula \cite{5, 
13, 15}. $H$ then has as its spectrum an infinite set of 
bands:
$\spec(H)=[E_0,E_1]\cup[E_2,E_3]\cup\cdots$. Let 
$\{\mu_n(x)\}_{n=1}^\infty$
be the eigenvalues of the Dirichlet Schr\"odinger operator 
in $L^2(x,x+L)$
(w.r.t.\ Lebesgue measure) with $u(x)=u(x+L)=0$ boundary 
conditions
$(E_{2n-1}\le \mu_n(x)\le E_{2n})$. The trace formula says 
that if $V$
is in $H^{1,2} ([0,L])$, where $H^{m,p}$ is the Sobolev 
space of distributions
with derivatives up to order $m$ in $L^p$, then
$$
V(x)=E_0+\sum_{n=1}^\infty (E_{2n}+
E_{2n-1}-2\mu_n(x)).\tag 2
$$

One of our main goals here is to prove a version of this 
trace formula
for arbitrary Schr\"odinger and Jacobi operators.

We will need the paired half-line Dirichlet operator 
$H\subd^x$ defined
on $L^2(-\infty,x)\oplus L^2(x,\infty)$ and $h\subd^n$ on
$l^2(\Bbb Z|m<n)\oplus l^2(\Bbb Z|m>n)$ with $u(x)$ (or 
$u(n)$\<)
vanishing boundary conditions. In the periodic case, it 
can be shown
that $\mu_n(x)$ are precisely the eigenvalues of 
$H\subd^x$ (as long as
$E_{2n-1}<\mu_n(x)<E_{2n}$, i.e., no equality).

The difference $(H-i)^{-1}-(H\subd^x-i)^{-1}$ is rank 1 
(and similarly
in the case of $h\subd^n$ if we define 
$(h\subd^n-i)^{-1}(n,m)\equiv 0)$\<)
and so trace class. As a result, the Krein spectral shift 
\cite{11}
exists; i.e., there is a function $\xi(x,\lambda)$ 
uniquely determined
a.e.\ in $\lambda$ w.r.t.\ Lebesgue measure by
$$\gather
\otr(f(H)-f(H\subd^x))=-\int_{-\infty}^\infty f\,'(\lambda)
\xi(x,\lambda)\,d\lambda,\tag 3\\
0\le\xi(x,\lambda)\le 1,\tag 4\\
\xi(x,\lambda)=0\quad\text{if}\ \lambda<\inf(\spec(H))
\endgather
$$
for any $C^1$ function, $f$, with $\sup_\lambda|(1+
\lambda^2)\,df/d\lambda|
<\infty$.

\<$\xi$ is a remarkable function which we claim is central 
to the proper
understanding of inverse problems; it will be discussed in 
detail in
three forthcoming papers which include detailed proofs of 
the theorems
that we present here \cite{6--8}. Our general trace 
formula is

\thm{Theorem 1S \cite6\rm} Let $V$ be continuous at $x$ 
and $E_0
\le\inf(\spec(H))$. Then
$$
V(x)=E_0+\lim_{\alpha\downarrow 0}\int_{E_0}^\infty 
e^{-\alpha\lambda}
(1-2\xi(x,\lambda))\,d\lambda.\tag 5S
$$
\ethm

\thm{Theorem 1J \cite6\rm} Let $E_-\le\inf(\spec(h))$ and 
$E_+
\ge \sup(\spec(h))$. Then
$$
v(n)=\frac 12 (E_-+E_+)+\int_{E_-}^{E_+}
\lf(\frac 12-\xi(n,\lambda)\rt)\,d\lambda.\tag 5J
$$
\ethm

\rem{Remarks} 1. If $V$ is smooth, there are higher-order 
trace 
relations including KdV invariants \cite7.

2. In the Jacobi case, $\xi(n,\lambda)=1$ if 
$\lambda>\sup(\spec(h))$,
which is needed for consistency in (5J). 

3. While we have singled out the Dirichlet boundary 
condition at 
$x\in\Bbb R$, any other selfadjoint boundary condition of 
the type
$\psi'(x)+\beta\psi(x)=0$, $\beta\in\Bbb R$, has been 
worked out
as well in \cite7.

4. Besides the motivating equation (2), two other special 
cases are in the
literature. Kotani and Krishna \cite{10} and Craig  \cite3 
discuss the case
where $V$ is bounded and continuous and (in our language) 
$\xi=\tfrac 12$ a.e.\
on $\spec(H)$; and Venakides \cite{16} has a trace formula 
when $V$ is positive
of compact support. In \cite6 we will discuss the relation 
of our work to
these in more detail. \endrem

\demo{Sketch of Proof} For simplicity, we consider only 
the Schr\"odinger
case and suppose $H\ge 0$ and take $E_0=0$. By (3)
$$
\otr(e^{-\alpha H}-e^{-\alpha H\subd^x})=\alpha
\int_0^\infty e^{-\alpha\lambda}\xi(x,\lambda)\,d\lambda.
$$
Moreover, a path integral argument shows that
$$
\otr(e^{-\alpha H}-e^{-\alpha H\subd^x})=\frac 12 
(1-\alpha V(x)+o(\alpha)).
$$
Given that
$$
\frac 12=\alpha\int_0^\infty e^{-\alpha\lambda}\frac 12 
\,d\lambda,\tag 6
$$
we get (5S) for $E_0=0$. \enddemo

A second critical result that we prove is

\thm{Theorem \rm2 \cite6\rm} For each $x\in\Bbb R$ and 
a.e.\ $\lambda$ in
$\Bbb R$, 
$$
\xi(x,\lambda)=\frac 1\pi\arg(G(x,x;\lambda+i0)).
$$
\ethm

\rem{Remark} $G$ is the integral kernel (resp.\ matrix 
elements)
of $(H-\lambda)^{-1}$ (resp.\ $(h-\lambda)^{-1}$\<). By 
general 
principles for each $x$, $\lim_{\varepsilon\downarrow 0} G
(x,x;\lambda+i\varepsilon)$ exists for a.e.\ $\lambda$. 
\endrem

\ex{Examples} 1. $V=0$. In the $H$ case, 
$G(x,x;\lambda)=(-\lambda)^{-1/2}$
for $\lambda\in\Bbb C\backslash[0,\infty)$ with the branch 
of square root,
so $G>0$ for $\lambda\in(-\infty,0)$. Thus, for 
$\lambda\in (0,\infty)$,
$G(x,x;\lambda+i0)=i|\lambda|^{-1/2}$ and 
$\xi(x,\lambda)\equiv \tfrac12$.
Equation (6) is then an expression of the known fact that $
\otr(e^{-\alpha H_0}-e^{-\alpha H_{\roman D,0}^x})=\tfrac 
12$ for all $\alpha$.

2. Let $V$ be periodic and in $H^{1,2}([0,L])$ with $V(x+
L)=V(x)$.
The spectrum of $H$ is $\bigcup_{n=0}^\infty[ E_{2n}, 
E_{2n+1}]$ as
noted already. Because $V$ is in $H^{1,2}([0,L])$,
$$
\sum_{n=0}^\infty|E_{2n}-E_{2n-1}|<\infty.\tag 7
$$
It can be shown (see, e.g., Kotani \cite9, Simon 
\cite{14}, and 
Deift and Simon \cite4) that $G(x,x;\lambda+i0)$ is pure 
imaginary
on $\spec(H)$, so $\xi=\tfrac 12$ there. Thus we claim 
(here and
below, we do not give a value to $\xi$ at points of 
discontinuity;
the real-valued function $\xi$ is only determined a.e.):
$$
\xi(x,\lambda)= \cases
\dfrac12, &\quad E_{2n}<\lambda<E_{2n+1},\\
1, &\quad E_{2n+1}<\lambda<\mu_{n+1}(x),\\
0, & \mu_{n+1}(x)<\lambda<E_{2n+2}, \endcases
$$
for $0\le\xi\le 1$, and $\xi$ jumps by $-1$ at $\mu_{n+
1}(x)$. 
Because of (7), 
$\int_{E_0}^\infty|1-2\xi(x,\lambda)|\,d\lambda<\infty$
and (5S) becomes (2).

3. Let $V(x)\to\infty$ as $|x|\to\infty$. Then $H$ has 
eigenvalues
$E_0<E_1<E_2<\cdots$ and $H\subd^x$ eigenvalues 
$\mu_1(x)<\mu_2(x)
<\cdots$ with $E_{n-1}\le\mu_n(x)\le E_n$. $|1-2\xi|=1$, 
so the 
integral in (5S) is not absolutely convergent if $\alpha$ 
is set
equal to zero and (5S) becomes a summability result; 
explicitly
$$
V(x)=E_0+\lim_{\alpha\downarrow 0}\alpha^{-1} 
\sum_{j=1}^\infty
[2e^{-\mu_j(x)\alpha} -e^{-E_j\alpha} -e^{-E_{j-1}\alpha}].
$$
For an explicit case, let $V(x)=x^2-1$ and place the 
Dirichlet condition
at $x=0$. Then
$$
E_n=2n,\qquad\mu_n(0)=
\cases 2n &\quad (n\ \roman{odd})\\
2(n-1) &\quad (n\ \text{even},\ n\ge 2),\endcases
$$
so $\xi(0,\lambda)=1$ on $(0,2)\cup (4,6)\cup\cdots$ and 
$\xi(0,\lambda)=0$
on $(2,4)\cup(6,8)\cup\cdots$ and formally
$$
\int_0^\infty(1-2\xi(0,\lambda))\,d\lambda=-2+2-2\cdots.
$$
The regularization (5S) is just the Abelian sum which is 
$-1$, 
which is exactly $V(0)$. 

4. Let $V(x)$ be short range in the sense that $V$ is 
$L^1(\Bbb R)$.
Then one can write down $\xi(x,\lambda)$ in terms of the 
reflection
coefficients $R(\lambda)$ and Jost functions $f_+
(x,\lambda)\ 
(\lim_{x\to\infty} e^{-i\lambda^{1/2} x} f_+
(x,\lambda)=1)$, viz \cite8
$$
\xi(x,\lambda)=\frac 12+\frac1\pi\arg \lf[\frac{1+
R(\lambda)f_+(x,\lambda)^2}
{|f_+(x,\lambda)|^2}\rt],\qquad \lambda>0.\tag 8
$$
In particular, $|\xi(x,\lambda)-\tfrac 12|\le\tfrac 12 
|R(\lambda)|$, 
and if $V\in H^{2,1}(\Bbb R)$, we have that
$$
\int_{E_0}^\infty\lf|\xi(x,\lambda)-\frac 
12\rt|\,d\lambda<\infty,\tag 9
$$
so
$$
V(x)=E_0+\int_{E_0}^\infty(1-2\xi(x,\lambda))\,d\lambda
$$
without a need for regularization.

5. There is a general summability result \cite8 like (9) 
also for the sum
of a smooth periodic potential and a sufficiently 
short-range potential
modeling impurity scattering in one-dimensional crystals. 
\endex

The Krein spectral shift has rather strong continuity 
properties:

\thm{Lemma 3a} Let $V_m(x)$ \RM(resp.\ $v_m(n)$\<\RM) 
converge to $V(x)$
uniformly for $x\in[-L,L]$ for each $L$ \RM(resp.\ to 
$v(n)$ for 
each $n$\<\RM) and so that $\inf_{x,m} V_m(x)<-\infty$ 
\RM(resp.\ 
$\sup_{n,m}|v_m(n)|<\infty$\<\RM). Then as measures in 
$\lambda$,
$\xi_m(x,\lambda)\,d\lambda$ converges weakly to 
$\xi(x,\lambda)\,
d\lambda$ for each fixed $x$. \ethm

It follows from Theorem 2 that

\thm{Lemma 3b} For each fixed $x$, 
$\spec_{ac}(H)=\{\lambda|0<\xi(\lambda,x)
<1\}^{-\roman{ess}}$ where $^{-\roman{ess}}$ is the 
essential closure.\ethm

Third, it follows from results of Kotani \cite9 in the 
Schr\"odinger
case and Simon \cite{14} in the Jacobi case:

\thm{Lemma 3c} If $V$ \RM(resp.\ $v$\RM) is periodic, then 
$\xi(x,\lambda)
\equiv\tfrac 12$ on $\spec(H)$ \RM(resp.\ 
$\spec(h)$\<\RM). \ethm

These three lemmas imply

\thm{Theorem 3 \cite6\rm} Suppose $V_m$ \RM(resp.\ 
$v_m$\<\RM) converge to
$V$ \RM(resp.\ $v$\<\RM) in the sense of Lemma \RM{3a} and 
each $V_m$
\RM(resp.\ $v_m$\<\RM) is periodic. Then for any 
measurable set
$S\subset\Bbb R$
$$
|S\cap\spec_{ac}(H)|\ge\varlimsup |S\cap \spec(H_m)|
$$
\RM(resp.\ replacing $H$ by $h$\<\RM) where $|\boldcdot|=$ 
Lebesgue measure.
\ethm

\ex{Example} Consider the Jacobi matrix with 
$v(n)=\lambda\cos(\pi\alpha n)$
(almost Mathieu or Harper's model). Avron et al.\ \cite1 
have proven
that if $\alpha$ is rational, then $|\spec(h_\alpha)|\ge 
4-2|\lambda|$.
Theorem 3 then implies (by approximating any $\alpha$ by 
rationals) that
$|\spec_{ac}(h_\alpha)|\ge 4-2|\lambda|$, slightly 
strengthening a recent result
of Last \cite{12}. In particular, we have a new proof of 
Last's 
spectacular result that 
$\spec_{ac}(h_\alpha)\ne\varnothing$ if 
$|\lambda|<2$ and $\alpha$ is a Liouville number. \endex

Finally, \cite6 will use $\xi$ to study the inverse problem.
Typical of our results is the following:

Let $V(x)\to\infty$ as $x\to\pm\infty$. Let $E_n(V)$ be the 
eigenvalues of $H=-d^2/dx^2+V$. We claim that when $V$ is 
even,
$\{E_n\}$ are a complete set of spectral data in the sense 
that

\thm{Theorem 4} If $V,W$ are continuous functions on $\Bbb 
R$ bounded
from below, going to infinity at $\pm\infty$, and obeying 
$V(x)=V
(-x)$ and $W(x)=W(-x)$ so that $E_n(V)=E_n(W)$ for all 
$n$, then
$V=W$. \ethm

Borg \cite2 proved this result over forty years ago. The 
$\xi$ 
function proof is natural, and we have an extension to
the nonsymmetric case. When $V$ is not symmetric, the 
Dirichlet
eigenvalues and the information about whether each is a 
Dirichlet
eigenvalue on $(-\infty,0)$ or $(0,\infty)$ also needs to 
be supplied.

\heading Acknowledgments\endheading
We thank H.~Kalf for valuable discussions. F.~Gesztesy is 
indebted
to the Department of Mathematical Sciences of the 
University of
Trondheim, Norway, and the Department of Mathematics at 
Caltech
for the hospitality extended to him in the summer of 1992.

\enddocument